\documentclass[11pt,a4paper,reqno]{amsart}
\usepackage{epsfig,amssymb,latexsym,amscd}

\usepackage[all,2cell]{xy}
\xyoption{v2}
\UseAllTwocells
\theoremstyle{plain}

\newtheorem{teo}{Theorem}[section]
\newtheorem{theo}[teo]{Theorem}
\newtheorem{coro}[teo]{Corollary}
\newtheorem{lema}[teo]{Lemma}
\newtheorem{prop}[teo]{Proposition}

\theoremstyle{remark}

\newtheorem{rk}[teo]{Remark}

\newtheorem{example}[teo]{Example}

\theoremstyle{definition}

\newtheorem{defi}[teo]{Definition}

\newcommand{\lcomod}{{}^C\!\mathcal M}
\newcommand{\rcomod}{\mathcal M^C}

\newcommand{\lmodA}{{}_{A}\mathcal M}
\newcommand{\rmodA}{\mathcal M_{A}}

\begin{document}
\title[Generator coalgebras are not necessarily quasi-coFrobenius]{Generator coalgebras are not necessarily quasi-coFrobenius}
\author{Mariana Haim}
\address{Facultad de Ciencias\\Universidad de la Rep\'ublica\\
Igu\'a 4225\\11400 Montevideo\\Uruguay\\}
\thanks{The first author would like to thank the Ministerio de Educaci\'on y Ciencia de Espa\~na}
\email{negra@cmat.edu.uy}
\author{Blas Torrecillas}
\address{Departamento de Algebra y An\'alisis Matem\'atico\\Universidad de Almer\'\i a\\
Cta Sacramento s/n, 04071 La Ca\~nada de san Urbano, Almer\'\i a \\Espa\~na\\}
\thanks{The second author would like to thank the projects $MTM2008-03339$.}
\email{btorreci@ual.es}

\begin{abstract} We study the problem of whether a coalgebra that generates its category of left (right) comodules is left (right) quasi-coFrobenius or not. We prove it does not hold in general, by giving a method of constructing counterexamples. This gives a negative answer to a question stated in \cite{kn:coalgen}. We also prove it is true for monomial pointed coalgebras and we characterize the quivers $Q$ such that $\Bbbk Q$ admits a monomial subcoalgebra that is left (right) quasi-coFrobenius.
\end{abstract}
\maketitle
\section{Introduction}\label{section:intro}

A coalgebra $C$ is said to be left (right) quasi-coFrobenius if it embeds, as a left (right) $C^*$-module, in a free left (right) $C^*$-module. It is well known that under this condition, the coalgebra generates the category of its left $C$- comodules. The question of whether the converse is true or not has been presented, in \cite{kn:coalgen}. It has been proved true under some finiteness conditions (for semiperfect coalgebras, for coalgebras whose coradical filtration is finite -see \cite{kn:coalgen}) and under some assumptions of duality (for coalgebras that are generators of their categories of left and right comodules, in particular for cocommutative coalgebras -see \cite{kn:rumanos}).\\

In this paper we prove that it does not hold in general, by constructing a family of counterexamples. We also prove that it holds for a class of coalgebras that includes the hereditary pointed coalgebras. More precisely, we prove it holds for subcoalgebras of path coalgebras that admit a (linear) basis of paths, the so called {\bf monomial pointed coalgebras}. We prove that a path coalgebra is quasi-coFrobenius if and only if its Gabriel quiver is discrete (i.e. with no arrows) and we give a list of all possible Gabriel quivers $Q$ such that $\Bbbk Q$ admits a subcoalgebra that is quasi-coFrobenius.\\

In what follows we present a brief description of the contents of this paper.\\
In Section \ref{section:preliminares} we recall a few basic notations and well known results in coalgebra
theory that will be used throughout the paper.\\
In Section \ref{section:relativo} we give the needed background to understand and treat the question. \\
In Section \ref{section:quivers} we present the results.
\section{Preliminaries and basic notations}\label{section:preliminares}

In this section we fix the basic notations and recall a few definitions and basic
concepts in coalgebra theory. We refer to \cite{kn:rumanos} and \cite{kn:sweedler} for more details.\\

The coalgebras we will work with will be $\Bbbk$-coalgebras, where $\Bbbk$ is a field.
For coalgebras and comodules we use Sweedler's notation. We denote by $\lcomod$ and $\rcomod$ the categories of left and right comodules over $C$ respectively. Similarly, if $A$ is an algebra over $k$, we denote by $\lmodA$ and $\rmodA$ the categories of left and right modules over $A$ respectively.\\
We denote by $\langle x \rangle ^C$ and ${}^C \langle x \rangle$ the respectively left and right comodules generated by the element $x\in C$ (i.e. the minimal left or right subcomodule of $C$ containing $x$).

\noindent
{\bf Finite localness of coalgebras}\\

Is is well known that simple right and left comodules are finite dimensional. As a consequence, every finitely generated right or left $C$-comodule is finite dimensional, and every simple coalgebra is finite dimensional. Moreover, every comodule is the union of its finite dimensional subcomodules.\\

\noindent
{\bf Injective comodules}\\

The category of right (left) comodules has enough injectives. Moreover, $C$ is injective as a right (left) $C$-comodule and for every right $C$-comodule $(M,\delta)$ there is an injection
$$
\delta: M \rightarrow M \otimes Coef(M) \cong Coef(M)^{dim M},
$$
where $Coef(M)$ is a subcoalgebra of $C$. If $M$ is simple, $Coef(M)$ is simple.\\

If $S$ is a simple right (left) comodule, we will call $E_r(S)$ (respectively $E_l(S)$) the injective envelope of $S$). When it is clear from the context, we just put $E(S)$.
Under this notation we have $E_r(S) \subseteq C$ for every right simple subcomodule $S$ of $C$ (the same for left simple comodules) and also
$$
C=\bigoplus_{S} E_r(S)=\bigoplus_{S} E_l(S),
$$
where $S$ belongs to the family of all simple right $C$-subcomodules of $C$ (respectively simple left $C$-comodules). It is important to remark here that every simple right $C$-comodule is isomorphic to a right $C$-subcomodule of $C$.\\



\noindent
{\bf Rational modules}\\

Every right (left) $C$-comodule can be thought as a left (right) $C^*$-module via the action $f\rightharpoonup m=\sum m_0 f(m_1)$ (repectively $m\leftharpoonup f= \sum f(m_{-1})m_0$).
In particular $C$ will be a right and a left $C^*$-module.\\
Left (right) $C^*$-modules that come in this way from right (left) $C$-comodules are called rational left (right) $C^*$-modules. The category of right (left) $C$-comodules is equivalent to the category of left (right) rational $C^*$-modules.

\section{Generator coalgebras and quasi-coFrobenius coalgebras}
\label{section:relativo}

We start by giving the needed background to understand the question proposed in \cite{kn:coalgen}, to which we will give a negative answer.\\

Let us understand first the notion of quasi-coFrobenius coalgebra. The following definition can be found in \cite{kn:rumanos}.
\begin{defi}
A coalgebra $C$ is said to be {\bf coFrobenius} if $C$ embeds in $C^*$ (as a left or as a right $C^*$-module).\\
A coalgebra $C$ is said to be {\bf left (right) quasi-coFrobenius} (shortly {\bf left (right) qcF}) if $C$, as a left (right)-$C^*$-module, can be embedded in a free left (right)-$C^*$-module.
\end{defi}

\begin{rk}\label{rk:qcF}
It is well known that the notion of coFrobenius is left-right symmetric (see for example \cite{kn:rumanos}).\\
On the other hand, the notion of quasi-coFrobenius  is not left-right symmetric (contrarely to the notion of quasi Frobenius ring). Actually, there are examples of coalgebras that are left qcF but not right qcF (see \cite{kn:rumanos}, example 3.3.7).
\end{rk}

Let us now go to the notion of generator.

\begin{defi}
Let ${\mathcal C}$ be a category with coproducts. An object $G$ in ${\mathcal C}$ is called a {\bf generator} of ${\mathcal C}$ if and only if for every object $X$ in ${\mathcal C}$ there is a set $I$ and an epimorphism $\varphi: G^{(I)} \rightarrow X$ in ${\mathcal C}$.
\end{defi}
\begin{defi}
We say that a coalgebra is  {\bf left f-qcF} if every finite dimensional (or, equivalently, finitely generated) rational left $C^*$-module (or equivalently, right $C$-comodule) embeds in a left free module.
\end{defi}

\begin{prop}\label{prop:relative}
\begin{itemize}
\item A coalgebra $C$ is left (right) qcF if and only if every left (right) rational $C^*$-module embeds in a free left (right) $C^*$-module.
\item A coalgebra $C$ is a generator of the category of its left (right) comodules if and only if it is left f-qcF.
\end{itemize}
\end{prop}
\begin{proof}
\begin{enumerate}
\item Take $M$ to be a left rational $C^*$-module, i.e. a right $C$-comodule. It is known that $M$ embeds in some $C^{(I)}$ as right $C$-comodule (and so as left $C^*$-module), for some $I$ and the direct implication follows. The converse is immediate.
\item Note that $M$ is a finitely generated left rational $C^*$-module if and only if it is a finitely generated right $C$-comodule, if and only if it is a finite dimensional right $C$-comodule. Then, see \cite{kn:coalgen}, Proposition 2.6.
\end{enumerate}
\end{proof}

\begin{rk}\label{rk:fgf}
From Proposition \ref{prop:relative} we conclude that every left (right) qcF coalgebra is a generator of the category of is left (right) $C$-comodules. The question we want to study is:
\begin{center}{\em Let $C$ be a coalgebra. If $C$ is a generator of $\lcomod$ then is $C$ left qcF?.}\end{center}
\end{rk}

Some positive results are known. For instance, if $C$ is a projective generator of $\lcomod$, then C is left and right qcF (\cite{kn:rumanos}, Corollary 3.3.11). In \cite{kn:qcfII}, it is shown that if $C$ is a generator for $\lcomod$ and for $\rcomod$ (in particular for cocommutative coalgebras), then C is left and right qcF. Moreover, in \cite{kn:coalgen}, the result is shown positive for coalgebras with finite coradical filtration and a characterization of qcF coalgebras is given. We give an alternative proof of this characterization in Proposition \ref{prop:semip}, in order to get closer to the problem.

\begin{defi}
A coalgebra $C$ is called {\bf left (right) semiperfect} if the injective envelope of every simple object in $\rcomod$ ($\lcomod$) is finite dimensional.
\end{defi}

\begin{prop}\label{prop:semip}
A coalgebra $C$ is left qcF if and only if it is a generator of $\lcomod$ and it is left semiperfect.
\end{prop}
\begin{proof}
For the direct implication we refer to \cite{kn:rumanos}, Corollary 3.3.6. For the converse, observe that if $S$ is a simple right $C$-comodule, then $E(S)$ is finite dimensional and we have that $E(S)$ embeds in a free left $C^*$-module (since $C$ is a generator in $\lcomod$ and therefore it is left f-qcF). Now, as $C=\bigoplus_S E(S)$, we conclude that $C$ embeds in some left free $C^*$-module and we are done.
\end{proof}

We finish this section with some general results that will be used in section \ref{section:quivers}.

\begin{prop}\label{prop:coideals}
A coalgebra $C$ is left f-qcF if and only if every finite dimensional right coideal in $C$ embeds (as a left $C^*$-module) in a free left $C^*$-module.
\end{prop}
\begin{proof}
The direct implication is immediate. For the converse, take $M$ to be a left finitely generated rational $C^*$-module. Then $M$ is finite dimensional and, as a right $C$-comodule, embeds in $M \otimes Coef(M) \cong Coef(M)^m$ where $m=dim_\Bbbk M$. As $Coef(M)$ is a finite dimensional subcoalgebra of $C$, it is a finite dimensional right coideal in $C$. Hence, $Coef(M)$ embeds in a free left $C^*$-module and so we conclude that $M$ embeds in a free left $C^*$-module.
\end{proof}

\begin{lema}\label{lemma:0}
Let $C$ be a $\Bbbk$-coalgebra and let $M$ be a right $C$-comodule. A map $F:M\rightarrow C^*$ is a morphism of left $C^*$-modules if and only if
$$
\sum F(m_0)(x) m_1= \sum x_1 F(m)(x_2), \forall m \in M, x \in C.
$$
\end{lema}
\begin{proof}
$F$ is a morphism of left $C^*$-modules if and only if $F(f\rightharpoonup m)=f.F(m), \forall m \in M, f \in C^*$. This is equivalent to
$$
F\left (\sum m_0 f(m_1)\right ) (x)= \left (f. F(m)\right )(x), \forall m \in M, x\in C,f \in C^*.
$$
So $F$ is a morphism or left $C^*$-modules if and only if $\sum F(m_0)(x) f(m_1)= \sum f(x_1)F(m)(x_2), \forall m\in M, x\in C, f \in C^*$. As this last equality is for every $f \in C^*$ we are done.
\end{proof}
\begin{coro}\label{coro}
 Let $M$ be a right $C$-comodule and $B$ be a $\Bbbk$-basis of $M$. Then $M$ is torsionless as a left $C^*$-module, if and only if for every $b\in B$ there exists a function $F_b:M \rightarrow C^*$ such that $F_b(b) \neq 0$ and $\sum F_b (m_0)(x) m_1= \sum x_1 F_b(m)(x_2), \forall m \in M, x \in C$.
\end{coro}
\begin{proof}
As $M$ is torsionless as a left $C^*$-module, there is a monomorphism of $C^*$ modules $F:M \rightarrow (C^*)^I$. If we denote by $F_i$ the composition of $F$ with the canonical projection to the $i$-th component in the diret product (for every $i\in I$), then each $F_i$ is also a morphism of left $C^*$-modules and therefore each $F_i$ verifies the formula presented in Lemma \ref{lemma:0}. \\
Take $b \in B$. As $b\neq 0$ we have that $F(b)\neq 0$ and therefore there is some $i\in I$ for which $F_i(b) \neq 0$. \\
For the converse it is enough to define $F:M \rightarrow (C^*)^{\# B}$ by putting, for each $b \in B$ as follows:\\ $F(b)$ is the vector in $(C^*)^{\# B}$ defined by $F_b(b)$ in the $b$-th coordinate and $0$ in all the other coordinates. \\
It is clear that $F$ is a morphism of $C^*$-modules that is injective.
\end{proof}

\section{Main results} \label{section:quivers}

We start by presenting path coalgebras and recalling an important result on pointed hereditary coalgebras. We refer to \cite{kn:chin}, \cite{kn:paths}, \cite{kn:pathspreprint} and \cite{kn:simson} for details.

\begin{defi}
Let $Q$ be a (oriented) quiver. In the vector space $\Bbbk Q$ of all linear combinations of paths in $Q$, a structure of coalgebra can be defined by:

$$
\begin{array}{ll}
\Delta: \Bbbk Q \rightarrow \Bbbk Q \otimes \Bbbk Q, & \Delta(p)=\sum_{xy=p} x\otimes y, \\
\ \\
\varepsilon:\Bbbk Q \rightarrow \Bbbk , & \varepsilon(p)=
\left \{
\begin{array}{ll}
1 & \mbox{ if $p$ is a trivial path}\\
0 & \mbox{ if not}
\end{array}
\right .
\end{array}
$$

The coalgebra defined in this way is called the {\bf path coalgebra} associated to $Q$ and will be denoted by $\Bbbk Q$.
\end{defi}

In the following $Q$ will be a quiver (with possible loops and multiple arrows), $C\subseteq \Bbbk Q$ will be a subcoalgebra of the path coalgebra and we will denote by $Q_0$ the vertices in $Q$ and by $Q_1$ the arrows in $Q$. If $p$ is a path in $\Bbbk Q$, we denote by $s(p)$ and $t(p)$ the vertices where $p$ starts and ends respectively. We also write sometimes $p: s(p)\rightarrow t(p)$ and call $s(p)$ and $t(p)$ the source and the target of $p$ respectively. Note that if $p$ is a path in $\Bbbk Q$ that belongs to $C$, then every arrow in $p$ belongs to $C$ (since $C$ is a subcoalgebra of $\Bbbk Q$). In particular $s(p), t(p) \in C$.

\begin{defi}
Let $\Bbbk Q$ be a path coalgebra and $C \subseteq \Bbbk Q$ a subcoalgebra. \\
We say that $C$ is {\bf admissible} if it contains all vertex and arrows in $Q$. \\
We say that an admissible subcoalgebra $C\Bbbk Q$  is {\bf monomial} of $\Bbbk Q$ if it is admissible and has a basis whose elements are all paths in $Q$.
\end{defi}

\begin{rk}
If $C$ is a subcoalgebra of a path coalgebra, there is only one quiver $Q$ such that $C$ is an admissible subcoalgebra of $\Bbbk Q$. This quiver is known as the {\bf Gabriel quiver} associated to $C$.
\end{rk}

\begin{example}\label{ex:quiver}
Consider the quiver
$
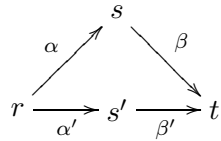
\xymatrix
{
& s \ar@{->}[rd]^{\beta} \\
r \ar@{->}[ru]^{\alpha} \ar@{->}[r]_{\alpha'}& s' \ar@{->}[r]_{\beta'} & t
}
$
The subspace of $\Bbbk Q$ with linear basis $\{r, \alpha, \alpha', s, s'\}$ is a monomial subcoalgebra of $\Bbbk Q$.\\
The subspace ot $\Bbbk Q$ with linear basis $\{\alpha \beta + \alpha'\beta', \alpha, \beta, \alpha', \beta', r, s, s', t\}$  is an admissible subcoalgebra of $\Bbbk Q $ that does not admit a basis formed by paths, so it is not monomial.
\end{example}

\begin{lema}\label{fuente}
Let $C$ be a monomial subcoalgebra of $\Bbbk Q$ and let $I\subseteq C$ be a left $C^*$-module that embeds in a free left $C^*$-module. Let $e\in I$ be a vertex. Then, there exists a path $p\in C$ such that:
\begin{itemize}
\item for all path $q \in I$ starting at $e$, $p=qp'$ for some path $p'\in C$, in other words: every path $q$ starting at $e$ belongs to the right coideal in $C$ generated by $p$,
\item for every non trivial path $r \in C$ ending at $e$, $rp \not \in C$.
\end{itemize}
\end{lema}
\begin{proof}
As $I$ is torsionless as a left $C^*$-module, there is a morphism $F:I \rightarrow C^*$ of left $C^*$-modules such that $F(e)\neq 0$. Let $B$ be a basis of $C$ whose elements are all paths and $p\in B$ be such that $F(e)(p)\neq 0$.\\
\begin{itemize}
\item Take $q\in I$ a path starting at $e$. Applying the equality of Lemma \ref{lemma:0} to $m=q, x=p$, we have
$$
\sum F(q_1)(p)q_2=\sum p_1 F(q)(p_2).
$$
Assume this is written as a linear combination of elements in $B$. Then, in the left term, the summand $F(e)(p)q$ is the only non zero multiple of $q$. Then the right term also has a summand that is a non zero multiple of $q$ and we are done.
\item Finally, assume that $r \in I$ is a non trivial path ending at $e$. Assume that $rp\in C$. We can then apply the equality to $m=e, x=rp$:
$$
F(e)(rp)e=\sum (rp)_1 F(e)\left ((rp)_2\right).
$$
But the right term has as non zero summand $r F(e)(p)$, while the left term is a multiple of $e \neq r$, arriving to a contradiction. Then $rp\not \in C$.
\end{itemize}

\end{proof}
\begin{rk}\label{rk:pintadep}
In the context of the previous lemma, we get in particular that $p$ verifies:
\begin{enumerate}
\item $s(p)=e$ (by taking $q=e$),
\item if there is a non trivial path $r\in C$ ending at $e$, then $p$ is also non trivial (otherwise $rp=r \in C$).
\end{enumerate}
\end{rk}

\begin{rk}\label{rk:quiver}
Let $C\subseteq \Bbbk Q$ be monomial. We recall that injective envelopes can be seen in terms of paths. Indeed, if $C\subseteq \Bbbk Q$ and $e\in C$ is a vertex the injective envelope in $C$ of the simple right (left)-comodule $\Bbbk e$ is the span of all paths in $C$ starting (ending) at $e$.
\end{rk}

We prove now that the implication holds in the case of monomial coalgebras.\\

\begin{theo}\label{theo:monomial}
Let $C\subseteq \Bbbk Q$ be a monomial subcoalgebra. If $C$ is left f-qcF, then $C$ is left qcF.
\end{theo}
\begin{proof}
In view of Proposition \ref{prop:semip}, it is enough to check that $C$ is left semiperfect. Let $Q_0$ be the set of vertices of the quiver $Q$. We will see that for every $e \in C\cap Q_0$, there are finitely many paths in $C$ starting at $e$ (see Remark \ref{rk:quiver}).\\
We start by proving that there is at most one arrow in $C$ starting at $e$. \\
Suppose we have two arrows $\alpha$ and $\beta$ starting at $e$. Take the vector space $M$ with basis $\{e,\alpha,\beta\}$. Then $M$ is a right coideal. As it is finite dimensional, it embeds in a free left $C^*$-module. By Lemma \ref{fuente}, there is some path $p\in C$ and some paths $a,b \in C$ such that $p=\alpha a=\beta b$. Then $\alpha=\beta$.\\
Suppose now that $C$ is not left semiperfect. Then some $E(\Bbbk e)$ is infinite dimensional. In particular, there is an arrow $\alpha$ starting at $e$.\\
Let $e'=t(\alpha)$. By Corollary $\ref{coro}$ applied to the right $C$-comodule $\Bbbk e'$, there exists a function $F_{e'}: \Bbbk e' \rightarrow C^*$ such that $F_{e'}(e') \neq 0$ and $F_{e'}(e')(x)e'=\sum x_1 F_{e'}(e')(x_2), \forall x \in C$.\\
By the proof of Lemma \ref{fuente} and by Remark \ref{rk:pintadep}, there is a path $x\in C$ such that $F_{e'}(e')(x)\neq 0$ and $s(x)=e'$. So $\alpha x$ is a path in $\Bbbk Q$ starting at $e$. As $C$ is not left semiperfect, there is a path in $C$ starting at $e$ that is longer than $\alpha x$, then it has to extend $\alpha x$ (since there is at most one arrow starting at each vertex). In particular, $\alpha x$ belongs to $C$.\\
To finish, observe that in the equality
$$
F_{e'}(e')(\alpha x)e'= \sum ( \alpha x)_1 F_{e'}(e')(\alpha x)_2,
$$
the left term is a multiple of $e'$ while the right term has a non zero summand (explicitly $\alpha F_{e'}(e')(x)$ ) which is not a multiple of $e'$, arriving to a contradiction.
\end{proof}

The following theorem gives in particular a characterization of all possible Gabriel quivers associated to monomial left(right) qcF coalgebras. We need one more technical result first.

\begin{lema}\label{lema:tecnico0}
Let $C\subseteq \Bbbk Q$ be a monomial subcoalgebra. Every finitely generated right coideal in $C$ embeds in a right coideal admitting a (linear) finite basis whose elements are all paths in $C$.
\end{lema}
\begin{proof}
Let $I$ be a right coideal. It is enough to take a generator $G$ of $I$ and write every element $g\in G$ as a linear combination of paths in $C$. If we call $G'$ the set of all paths obtained for all $g\in G$ in this way, we get that the ideal $I'$ generated by $G'$ contains $I$. Taking all the paths generated by $G'$ we obtain a linear generator of $I'$ than can be reduced to a linear basis of $G'$ (whose elements are all paths). It is clear that if $G$ is finite, the obtained basis is also finite.
\end{proof}

\begin{theo}\label{theo:gabriel}
Let $C \subseteq \Bbbk Q$ be a monomial subcoalgebra. Then $C$ is left qcF if and only the following two conditions hold:
\begin{enumerate}
\item $Q$ is a disjoint union of quivers of the following type:
\begin{itemize}
\item trivial quivers (a single point and no arrows),
\item quivers such that there is exactly one arrow starting at each vertex,
\end{itemize}
\item for each vertex $e \in Q_0$, there are finitely many paths starting at $e$ that belong to $C$.
\end{enumerate}
\end{theo}
\begin{proof}
Assume first that $C$ is left qcF. Then, from the proof of Theorem \ref{theo:monomial} we know that there is at most one arrow starting at each vertex. \\
Let $e$ be a vertex such that no arrow starts at $e$. By Remark \ref{rk:pintadep}, (2), there are no arrows ending at $e$ and so the connected component of $e$ is the trivial quiver.\\
Take a non trivial connected component of $Q$. There is exactly one arrow starting at every vertex and we are done.\\
As $C$ is left qcF, the injective envelope $E(\Bbbk e)$ is finite dimensional, so there are finitely many paths starting at each vertex (see Proposition \ref{prop:semip}).\\
For the converse, we can assume, without loss of generality, that $Q$ is connected (since a direct sum of left qcF coalgebras is left qcF).\\
Moreover, in view or Theorem \ref{theo:monomial}, it is enough to prove that if $Q$ is the two conditions are verified, then $C$ is left f-qcF. \\
If $Q$ is the trivial quiver, then $C\cong k$ and it is trivially coFrobenius.\\
For the other case, we will prove that any finite dimensional right coideal embeds, as a left $C^*$-module, in some free left module. Let $I\subseteq C$ be a finite dimensional right coideal. We know by Lemma \ref{lema:tecnico0} that $I$ embeds in a right coideal $I'$ that admits a linear basis $B=\{p_1,p_2,\cdots,p_n\}$ whose elements are all paths. Then we can put $B=B_1\cup B_2 \cup \cdots \cup B_k$, where for each $i$, the paths in $B_i$ share the source (that we call $s_i$) and for each pair $i \neq j$, $s_i \neq s_j$. In this situation, we get $I=I_1\oplus I_2 \oplus \cdots \oplus I_k$, where each $I_i$ is linearly generated by $B_i$. It is easy to see that each $I_i$ is a (finite dimensional) right coideal.\\
We can then assume that we have a finite dimensional right coideal whose elements share the source $s$. Take $p$ to be the path starting at $s$ of maximal length and let $t=t(p)$. Take a basis $B_I$ of $I$ and for each path $b\in B_I$, define $b'$ to be the (only) path from $t(b)$ to $t$. Let $B'_I$ be the set of all $b'$. Clearly $B'_I$ is linearly independent, so extend it to a basis $B$ of $C$ and take $B^*$ to be its dual basis.\\
Define $F:I \rightarrow C^*$ by $F(b)=(b')^*$. Clearly $F$ is injective. It is easy to verify that the formula of Lemma \ref{lemma:0} holds for $F$. \\
\end{proof}
\begin{rk}
Of course, Theorem \ref{theo:gabriel} admits a symmetric version for right qcF monomial coalgebras.
\end{rk}

\begin{example}
Take the quiver $A_\infty$ and call $x_1,x_2,\cdots,x_n,\cdots$ its (ordered) vertices. Define $C\subseteq \Bbbk A_\infty$ to be the subcoalgebra formed by all paths starting at $x_n$, of length less or equal to $n$. \\
Then $C$ is left qcF. Note that $C$ has no finite coradical series (since there are arbitrary long paths in $C$ -see Remark \ref{rk:quiver}).
\end{example}

Recall that a coalgebra $C$ is said to be
\begin{itemize}
\item {\bf pointed}, if all its simple subcoalgebras are one dimensional,
\item {\bf hereditary}, if the quotient of homomorphic images of injective right (or left) comodules are injective.
\end{itemize}
The following result is well-known and very useful in order to study pointed hereditary coalgebras.

\begin{theo}\label{theo:pointed}
If $Q$ is a quiver, then every subcoalgebra of $\Bbbk Q$ is a pointed coalgebra, (whose simple subcoalgebras are the generated by the vertices).\\
Moreover, $\Bbbk Q$ is pointed and hereditary.  Conversely, any pointed and hereditary coalgebra $C$ is isomorphic to a certain path coalgebra $\Bbbk Q$. More precisely, there is a unique quiver $Q$ (called the Gabriel quiver of $C$) such that $C$ is isomorphic to $\Bbbk Q$.
\end{theo}
\begin{proof}
See \cite{kn:chin}.
\end{proof}

\begin{rk}
It could be relevant to remark at this stage that in the case that the base field $\Bbbk$ is algebraically closed, any coalgebra $C$ is equivalent (in the Morita sense, i.e. the comodules categories will be equivalent as categories) to a subcoalgebra of a pointed coalgebra (see \cite{kn:chinmont}).
\end{rk}

\begin{coro}
Let $C$ be an hereditary pointed coalgebra with Gabriel quiver $Q$. The following are equivalent:
\begin{enumerate}
\item $C$ is left f-qcF,
\item $C$ is right f-qcF,
\item $C$ is left qcF,
\item $C$ is right qcF,
\item $C$ is coFrobenius,
\item $C$ is cosemisimple,
\item $Q$ is dicrete.
\end{enumerate}
\end{coro}
\begin{proof}
We have proved in Theorem \ref{theo:monomial} that (1) is equivalent to (3). \\
Assume (3). As $C=\Bbbk Q$, we have that all paths in $Q$ are elements in $C$. If $C$ is left qcF, then it is left semiperfect, so there are finitely many paths in $\Bbbk Q$ starting at each vertex. Using Theorem \ref{theo:gabriel}, we get that the only possibility for $Q$ is to be a disjoint union of trivial quivers, i.e. a discrete quiver. So we have that (3) implies (7).
Now, if $Q$ is discrete and $C=\Bbbk Q$, then the map $\varphi:C \rightarrow C^*$, defined by $\varphi (e)(f)=\delta_{e,f}$ is an injective morphism of left (and right) $C^*$-modules, so $C$ is coFrobenius and we have (7) implies (5).
Clearly (5) impies (1), so we have that (1), (3), (5) and (7) are equivalent. Similarly, (2), (4) , (5) and (7) are equivalent.
Finally, it is known that (6) implies (5) (see for example \cite{kn:rumanos}) and it is easy to see that (7) implies (6) (since $C=\Bbbk Q$ will be the direct sum of the -simple- coalgebras formed by the vertices in $Q$).
\end{proof}

Observe that, in example \ref{ex:quiver}, the basic elements which are not paths keep the property of being linear combinations of paths sharing the source and the target. The next Proposition asserts that, in general, every subcoalgebra of a path coalgebra has this property.

\begin{defi}
We will say that a linear combination of paths is a right (left) {\bf multipath} if the paths in it share the target (the source). If $p$ is a right (left) multipath, we denote by $t(p)$ ($s(p)$) these common target (source) and we say that $p$ ends at $t(p)$ (starts at $s(p)$). A right (left) multipath $q$ will said to be trivial if it is a scalar multiple of a vertex.\\
We say that $p$ is a multipath if it is a left and right multipath.
\end{defi}

\begin{rk}\label{rk:quiver2}
Let $C\subseteq \Bbbk Q$ be a subcoalgebra and $e \in C$ be a vertex. The injective envelope of $\Bbbk e$ as a right (left) simple comodule is spanned by all multipaths in $C$ starting (ending) at $e$.
\end{rk}

\begin{prop}\label{prop:multipaths}
Let $Q$ be a quiver. If $C \subseteq  \Bbbk Q$ is a subcoalgebra, there exists a linear basis of $C$ whose elements are multipaths.\\
Moreover if $I \subseteq \Bbbk Q$ is a right (left) coideal, there exists a linear basis of $I$ whose elements are right (left) multipaths.
\end{prop}
\begin{proof}
See \cite{kn:pathspreprint}, Theorem 2.4 for the case of subcoalgebras and note that the proof can be adapted to have the result on right or left coideals.
\end{proof}

The following definition, lemma and corollary are a tool that will allow us to restrict our attemption to right coideals that admit a (finite) basis of multipaths. The notion of {\bf independent} set is introduced in order to obtain "smaller" counterexamples but it is not essential to solve the problem.
\begin{defi}
Let $C\subseteq \Bbbk Q$ be an admissible subcoalgebra. Call $M(C)$ the set of all multipaths in $C$. A subset $F \subseteq M(C)$ will be said to be {\bf independent} if it is linearly independent and for each $x\in F$,
$$
\langle x \rangle ^C \cap F =\{x\}.
$$
\end{defi}

\begin{lema}\label{lema:tecnico}
Let $C\subseteq \Bbbk Q$ be a coalgebra.
\begin{enumerate}
\item Every finite dimensional right coideal in $C$ embeds in a right coideal generated by multipaths.
\item Let $B$ be a generator of a finite generated right coideal $I$, whose elements are multipaths in $C$. Then there is a finite independent set $F\subseteq B$ that generates $I$.
\end{enumerate}
\end{lema}
\begin{proof}
\begin{enumerate}

\item Let $I$ be a finite dimensional right coideal in $C$ and $G$ be a generator of $I$ whose elements are right multipaths (see Proposition \ref{prop:multipaths}). Each right multipath $g \in G$ is of the form $\alpha^g_1 + \alpha^g_2 \cdots \alpha^g_{n_g}$, where each $\alpha^g_i$ is a multipath in $C$. Take the right coideal $I'$ generated by $\{\alpha^g_i \mid g\in G, 1\leq i\leq n_g\}$. Then $I \subseteq I'\subseteq C$ and we are done

\item Let $B=\{x_1,x_2,\cdots,x_n\}$. We start by removing from $B$ all the elements in $\langle x_1 \rangle \cap B$ that are different from $x_1$, obtaining a new set $B_1=\{x_1,x^1_2, x^1_3, \cdots, x^1_{n_1}\}\subseteq B$ that keeps generating $I$. Then, we remove from $B_1$ all the elements in $\langle x^1_2 \rangle \cap B_1$ different from $x^1_2$ and obtain a new set $B_2=\{x^1_2, x^2_2, \cdots, x^2_{n_2}\}\subseteq B_1$ that keeps generating $I$. We continue in the same way, by considering in $B_2$ a new element (different from $x_1$ and $x^1_2$) and obtaining a $B_3\subseteq B_2$ that generates $I$. As $B$ is finite and the cardinal of the generators decreases, the process stops and it is easy to see that if we reduce the last set to a linearly independent one, the obtained set is finite and independent.
\end{enumerate}
\end{proof}

\begin{coro}\label{cor:tecnico}
Let $I$ be a finitely generated right coideal in $C\subseteq \Bbbk Q$. There exists a right coideal $I'\subseteq C$ containing $I$ and generated by a finite independent set of multipaths in $C$.
\end{coro}
\begin{proof}
It follows immediately from Lemma \ref{lema:tecnico}.
\end{proof}
We describe now the construction that will allow us to give examples of coalgebras that do not verify the implication.

Let $Q$ be a quiver and $C\subseteq \Bbbk Q$ be a coalgebra. \\

Extend $Q$ to a quiver $\overline{Q}$ by adding, for each finite independent set of multipaths $F=\{p_1,p_2,\cdots p_n\} \subseteq M(C)$:
\begin{itemize}
\item a vertex that we call $e_F$,
\item $n$ different arrows $\alpha_{1,F}, \alpha_{2,F},\cdots,\alpha_{n,F}$, where for each $i\in \{1,2,\cdots n\}$, $s(\alpha_{i,F})=t(p_i), t(\alpha_{i,F})=e_F$,
\end{itemize}
 i.e. the set of vertices is $\overline{Q}_0=Q_0=Q \cup \{e_F\mid F \subseteq M(C) \mbox{ finite and independent}\}$
 and the set of arrows is $\overline{Q}_1=Q_1 \cup \left \{\alpha_{i,F}\mid F \subseteq M(C) \mbox{ finite}, i \in \{1,2,\cdots,\# F\}\right \}$.\\
For each finite $F=\{p_1,p_2,\cdots, p_n\}\subseteq M(C)$, define $q_F$ to be the multipath $p_1\alpha_{1,F}+p_2\alpha_{2,F}+\cdots + p_n\alpha_{n,F}\in \Bbbk \overline Q$.\\
Then define the coalgebra $\overline C$ as the subcoalgebra of $\Bbbk \overline Q$ generated by the set $C\cup \{ q_F \mid F=\{p_1,p_2,\cdots, p_n\}\subseteq M(C) \mbox{ finite and independent}, n\in \mathbb N\}$.

Now, for a subcoalgebra $C\subseteq \Bbbk Q$, we define, for each $n\in \mathbb N$, a coalgebra $T^n(C)$, by recursion, as follows:
$$
T^0(C)=C, \hspace{3mm}, T^{n+1}(C)=\overline {T^n(C)}.
$$
Note that $T^0(C)\subseteq T^1(C)\subseteq T^2(C)\subseteq \cdots \subseteq T^n(C)\subseteq \cdots$.

\begin{defi}
Let $C \subseteq \Bbbk Q$ be a subcoalgebra. Under the above notation, the coalgebra $T(C)=\bigcup_n T^n(C)$ will be called the {\bf tail closure of $C$}.
\end{defi}

\begin{prop}\label{prop:T(C)semip}
If $C$ is not left semiperfect, then $T(C)$ is not left semiperfect.
\end{prop}
\begin{proof}
Take $e \in Q$ to be a vertex such that there are infinitely many multipaths in $C$ starting at $e$ (note that such a vertex exists because $C$ is not left semiperfect -see Remark \ref{rk:quiver2}). As $C\subseteq T(C)$, there are infinitely many paths in $T(C)$ starting at $e$ and therefore $T(C)$ is not left semiperfect.
\end{proof}

The rest of the work is devoted to prove that $T(C)$ is a generator of the category of its left comodules. \\

\begin{rk}\label{rk:circuits}
Call $T(Q)=\bigcup_n T^n(Q)$, where $T^0(Q)=Q$ and $T^{n+1}(Q)=\overline{T^n(Q)}$. Then it is clear that $T(C) \subseteq \Bbbk T(Q)$.
\end{rk}

\begin{prop}\label{prop:cyclic}
Every finite dimensional right coideal in $T(C)$ embeds in a right coideal generated by an element of the form $q_F$.
\end{prop}
\begin{proof}
Let $I$ be a finite dimensional right coideal in $T(C)$. From Corollary \ref{cor:tecnico}, there exists a right coideal $I'$ in $T(C)$ such that $I\subseteq I´$ and $I'$ is generated by a finite independent set $F$ of multipaths. Moreover, by the construction of $F$, we know that for some $n\in \mathbb N$, $F\subseteq T^n(C)$ and so $q_F$ (as constructed before) is an element in $T^{n+1}(C)\subseteq T(C)$. By the construction of $q_F$, it is clear that  $F\subseteq \langle q_F\rangle ^C$ (the right coideal generated by $q_F$). So $I \subseteq I' \subseteq \langle q_F \rangle ^C$.
\end{proof}

The following lemma admits a more general version on rings, but we state here only the result we will need later.

\begin{lema}\label{lemma:anulador}
Let $q\in C$ and $I=\{f \in C^* \mid f\rightharpoonup q=0\} \subseteq C^*$. If there is some $X\subseteq C^*$ such that $I=l(X)$, then the cyclic left $C^*$ module generated by $q$ embeds in a free left $C^*$-module.
\end{lema}
\begin{proof}
Take the morphism of left $C^*$-modules $\varphi:C^* \rightarrow (C^*)^{\# X}$ defined by $\varphi(\varepsilon)=(x)_{x\in X}$ (recall that $\varepsilon$ is the identity element in $C^*$). It is clear that $ker (\varphi)=l(X)=I$, so $I$ is a left ideal in $C^*$ and $\frac{C^*}{I}$ embeds in $(C^*)^{\# X}$. \\
Note that $\langle q \rangle ^C \cong \frac{C^*}{I}$ (as left $C^*$-modules), so $\langle q \rangle ^C$ embeds in $(C^*)^{\# X}$ as a left $C^*$-module. Now, $\langle q \rangle ^C$ is a finite dimensional rational left module, so it is a finitely cogenerated left $C^*$-module, hence for some $n\in \mathbb N$, we get that $\langle q \rangle ^C$ embeds in $(C^*)^n$.
\end{proof}

\begin{theo}\label{theo:salio}
Let $Q$ be a quiver and $C\subseteq \Bbbk Q$ be an admissible subcoalgebra. Then its tail closure $T(C)$ is a generator of the category of its left comodules.
\end{theo}
\begin{proof}
In view of propositions \ref{prop:coideals} and \ref{prop:cyclic}, it is enough to prove that every right coideal generated by an element of the form $q_F$ embeds in some free left $C^*$-module. \\
Let $I_F=\{f \in T(C)^* \mid f \rightharpoonup q_F=0\}$ and call (as before) $e_F=t(q_F)$.\\
Take $X=\{f \in T(C)^* \mid f(p)=0 \mbox{ for every multipath } p \mbox{ such that } s(p)\neq e_F\}$. We claim that $I_F=l(X)$.\\
Note first that $I_F=\{f \in C^* \mid \mbox{ such that }f(x)=0, \mbox{ for all multipath } x \in \ ^C\langle q_F \rangle\}$. By the construction of $T(C)$, the reader can notice that the multipaths in $^C\langle q_F \rangle $ are exactly the multipaths ending at $e_F$.\\
Now, assume that $f \in I_F$ and $g \in X$. Take $r$ to be a multipath in $T(C)$. We want to prove that $(fg)(r)=0$. But
in $(fg)(r)=\sum f(r_1)g(r_2)$, we have $f(r_1)=0$ for each summand such that $r_2$ starts in $e_F$ (since in this case $r_1$ ends in $e_F$) and $g(r_2)=0$ for each summand such that $r_2$ does not start in $e_F$. Then $(fg)(r)=0$.\\
Conversely, suppose that $(fg)(r)=0, \forall g \in X, \forall r\in C$. We want to prove that $f(x)=0$ for all multipaths $x$ ending at $e_F$. Take such an $x$ and put $G=\{x\}$. Take $q_G$ as defined in the construction of $T(C)$, so $q_G=x \alpha_G$, where $\alpha_G \in T(C)$ is an arrow starting at $e_F$ and ending in $e_G$. Now, $(fg)(q_G)=0$, but the only finally non zero summand in $(fg)(q_G)$ is $f(x)g(\alpha_G)$, so $f(x)g(\alpha_G)=0, \forall g \in X$. As $s(\alpha_G) = e_F$, there is some $g\in X$ such that $g(\alpha_G) \neq 0$. Hence, $f(x)=0$.\\
The proof finishes by using Lemma \ref{lemma:anulador} to get that $\langle q_F\rangle ^C$ embeds in some free left $C^*$-module.
\end{proof}

\begin{coro}\label{coro:finally!}
If $C$ is a coalgebra that is not left semiperfect, then $T(C)$ is a coalgebra that generates the category of its left comodules and is not semiperfect.
\end{coro}
\begin{proof}
It follows from Proposition \ref{prop:T(C)semip} and Theorem \ref{theo:salio}.
\end{proof}
\begin{rk}
Note that the coalgebra $\Bbbk A_\infty$ is not semiperfect, so $T(\Bbbk A_\infty)$ is a particular counterexample for the question. (This is in some sense the {\em simplest} counterexample this construction can give.)
\end{rk}
\begin{rk}
It is interesting to note that this question, stated in the language of comodules (or of rational modules), can be stated also for all $C^*$-modules, obtaining the following:
{\em Let $C$ be a coalgebra such that every finitely generated left $C^*$-module embeds in a free left $C^*$-module. Is it true that then every left $C^*$-module embeds in a free left $C^*$-module?}\\
This question is known in ring theory as the $FGF$-conjecture (and formulated as {\em Is every left $FGF$ ring a $qF$ ring?}). \\
In this sense, the problem we have treated here, is a relative version to rational modules of the $FGF$-conjecture for modules over pseudocompact algebras.
\end{rk}


\begin{thebibliography}{99}

\bibitem{kn:chin} Chin, W., {\em Hereditary and path coalgebras}, Communications in Algebra, Vol. 30 (4), pp 1829-1831, (2002).

\bibitem {kn:chinkl} Chin, W., Kleiner, M., Quinn, D., {\em Almost split sequences for comodules} Jounal of Algebra, Vol 249 (1), pp 1-19, (2002).

\bibitem{kn:chinmont} Chin, W, Montgomery, S., {\em Basic coalgebras.} Modular interfaces (Riverside, CA, 1995), AMS/IP Studies in Advanced Mathematics, Vol 4, pp 41-47 (1997).


\bibitem{kn:rumanos} Dascalescu, S., Nastasescu, C., Raianu, S., {\em Hopf algebras. An Introduction}, Pure and applied mathematics, A Series of Monographs and Textbooks, New York, (2000).


\bibitem{kn:qcfII} G\'omez Torrecillas, J., Manu, C., Nastasescu, C., {\em Quasi-co-Frobenius Coalgebras II}, Communications in Algebra, Vol. 31 (10), pp 5169-5177, (2003).

\bibitem{kn:paths} Jara, P., Merino, L., Navarro, G. {\em On path coalgebras of quivers with relations}, Colloquium Mathematicum, Vol. 102 (1), pp 49-65, (2005).

\bibitem{kn:pathspreprint} Jara, P., Merino, L. Navarro, G., Ruiz, J.F., {\em Prime path coalgebras}, preprint, (2004).


\bibitem{kn:coalgen} Nastasescu, C., Torrecillas, B. Van Oystaeyen, F. {\em When is a coalgebra a generator?}, Alg. Representation theory, Vol.11 (2), pp 179--190 (2008).

\bibitem{kn:nicholson} Nicholson, W.K., Youssif, M;.F. {\em Quasi-Frobenius rings}, Cambridge Tracts in Mathematics, Vol. 158, Cambridge University Press (2003).


\bibitem{kn:simson} Simson, D., {\em Coalgebras, comodules, pseudocompact algebras and tame comodule type}, Colloquium Mathematicum, Vol. 90 (1), pp 101-150, (2001).

\bibitem{kn:sweedler} Sweedler, M. {\em Hopf algebras}, New York: W.A. Benjamin, Inc. 336 p (1969).

\end{thebibliography}
\end{document}